\documentclass[leqno]{amsart}

\usepackage{amsfonts}
\usepackage{amssymb, latexsym, amsmath, pb-diagram}
\usepackage{pstricks-add}
\usepackage{lamsarrow, pb-lams}
\usepackage{multicol}
\usepackage{bm}
\usepackage[mathscr]{eucal}
\usepackage{xy}
\usepackage{epic,eepic}
\xyoption{all}

\numberwithin{equation}{section}

\def\w{\widetilde}
\def\wh{\widehat}

\def\nr{\refstepcounter{thm}\thethm}
\def\cf{{\it cf.}\ }
\def\o{\overline}

\newcommand{\nn}{\mathbb{N}}

\newcommand{\PP}{\mathbb{P}}
\newcommand{\QQ}{\mathbb{Q}}

\newcommand{\kk}{\mathbf{k}}
\newcommand{\ba}{\mathbf{a}}
\newcommand{\xx}{\mathbf{x}}

\newcommand{\ZZ}{\mathcal{Z}}
\newcommand{\II}{\mathcal{I}}

\begin{document}
\title[Simplicial complements and  Moment-angle Complexes]
{\bf The homology of simplicial complements and the cohomology of moment-angle complexes}
\author[Q. Zheng \& X. Wang]{ Qibing Zheng and Xiangjun Wang}
\thanks{The authors were supported by NSFC grant No. 10771105 and No. 11071125}
\address{School of Mathematical Sciences and LPMC, Nankai University,
Tianjin, 300071, P.R.China.}
\email{xjwang@nankai.edu.cn}
\address{School of Mathematical Science and LPMC, Nankai University,
Tianjin, 300071, P.R.China}
\email{zhengqb@nankai.edu.cn}
\maketitle

\begin{abstract}
In this paper, we use the Taylor resolution of a monomial ideal $I$ to compute the algebra $Tor_*^{\kk[\xx]}(I,\kk)$. As an application, we define the simplicial complement $\PP$ of a simplicial complex $K$ and construct a chain complex $(\o\Lambda(\PP), d)$ to compute the algebra $Tor_*^{\kk[\xx]}(\kk(K),\kk)$. We also give a description of the cohomology of generalized moment-angle complex with respect to the homology of the chain complex $H_*(\o\Lambda(\PP),d)$.
\end{abstract}

\section{Introduction and statement of results}

   The moment-angle complexes have been studied by topologists for many years
(\cf \cite{Porter} \cite{Lopez}). In 1990's
Davis and Januszkiewicz \cite{DJ} introduced toric manifolds theory which is studied
intensively by algebraic geometers. They observed
that every quasi-toric manifold is the quotient of a moment-angle complex by the free action
of a real torus, here the moment-angle complex is denoted by $\mathcal{Z}_K$ corresponding to
an abstract simplicial complex $K$.
The topology of $\ZZ_K$ is complicated and getting more attentions by topologists lately
(\cf \cite{GM} \cite{Hochster} \cite{Baskakov} \cite{Panov} \cite{Franz}).
Recently a lot of work has been done on generalizing the moment-angle complex
$\ZZ_K=\ZZ_K({D^2}, {S^1})$ to pairs of spaces $({X}, {A})$
(\cf  \cite{BB}, \cite{BBCG08}, \cite{GTh}, \cite{LP}).
In this paper we study the cohomology of the generalized moment-angle complexes
$\ZZ_K({X}, {A})$ corresponding to the pairs of spaces
$({X}, {A})$ with inclusions $A_i\hookrightarrow X_i$ being
homotopic to constant for all $i$.

The homology algebra $Tor_{*}^{\kk[\xx]}(\kk(K),\kk)$ of the Stanley-Reisner {\it face ring} $\kk(K)$ of the simplicial complex $K$
is the cohomology of $\ZZ_K$ (\cf \cite{BP}, \cite{Stanley}). The usual way of computing this
$Tor$ group is Hochster formula. But Hochster formula gives no information on the product of the algebra. The reason is that the usual Hochster formula is obtained from the tensor product of $\kk(K)$ with the Koszul complex, which is the minimal resolution of $\kk$. But this chain complex still a module over $\kk(K)$, the actual computation is quite complicated. In this paper, we use Taylor resolution (Definition 2.5) of the simplicial complement $\PP$ (Definition 2.4) of $K$ to get a DGA $(\o\Lambda(\PP),d)$ that give a new proof of the Hochster formula and determines the algebra structure of $Tor_{*}^{\kk[\xx]}(\kk(K),\kk)$ (Theorem 2.6). The advantage of the method is that the DGA algebra $(\o\Lambda(\PP),d)$ is a vector space over $\kk$ and the product of it directly determines the product of $Tor_{*}^{\kk[\xx]}(\kk(K),\kk)$. This result is the special case Theorem 2.2, the Taylor resolution of a monomial ideal.
   Let $\PP=\{\sigma_1, \sigma_2, \cdots , \sigma_s\}$  be a sequence of subsets of
$[m]=\{1,2,\cdots,m\}$, which is called a {\it simplicial complement} in this paper. It is easy to see that
$$\II_\PP=\langle \xx_{\sigma_1}, \xx_{\sigma_2}, \cdots, \xx_{\sigma_s}\rangle$$
is an ideal generated by square-free monomials. Thus there is an unique simplicial complex $K_\PP$
such that $\II_{K_\PP}=\II_\PP$.

Based on the computation of $Tor$ algebras and the decomposition of $\Sigma\ZZ_{K}({X}, {A})$ given by
Bahri, Bendersky, Cohen and Gitler (\cf \cite{BB,BBCG08}), we proved in Theorem 3.7 the cohomology of the generalized moment-angle complex $\Sigma\ZZ_{K}({X}, {A})$ which is related very much to the group $Tor_{*}^{\kk[\xx]}(\kk(K),\kk)$. As applications, we also consider the cohomology of $\ZZ_{K_\PP}({X}, {A})$
for some special triples of $CW$-complexes $({X}, {A})$ including
all the $X_i$ are contractible; all the $A_i$ are contractible and
$({X}, {A})=({S^2}, {S^1})$.

{\it Acknowledgements } The authors are grateful to Z.~L\"u for his introduction of this topic. The authors
are indebted to V.~M.~Buchstaber for his comments and suggestions, especially for the name of definition
{\it simplicial complement} which is originally called {\it partition.} The authors also thank S. Gitler for his
helpful suggestions.

\section{The homology of systems of generators and the $Tor$ algebra of face ring}

    In this paper, $[m]=\{1,2,\cdots, m\}$ and $2^{[m]}$ is the set of all subsets of $[m]$.
$\kk$ is a field. $\kk[\xx]$ is the polynomial ring  on $m$ indeterminates
$\xx=\{x_1, x_2, \cdots, x_m\}$. We always regard $\kk[\xx]$ as a $\nn^m$ graded algebra over $\kk$, where
$\nn=\{0,1,2,\cdots\}$. Precisely, $\kk[\xx]$ has a base of monomials expressed as
\[
\xx^\ba=x_1^{a_1}x_2^{a_2}\cdots x_m^{a_m}
\]
for every vector $\ba=(a_1,\cdots,a_m)\in \nn^m$. For vectors $\ba_i=(a_{i,1},\cdots,a_{i,m})\in\nn^m$,
$i=1,\cdots,n$, the least common multiple $\ba=[\ba_1,\cdots,\ba_n]$ is the vector $(a_1,\cdots,a_m)$ with
$a_j={\rm max}\{a_{1,j},\cdots,a_{n,j}\}$ for $j=1,\cdots,m$. For one monomial $\ba$, define $[\ba]=\ba$.
\vspace{2mm}

  {\bf Definition 2.1} {\it A sequence $\PP=\{\ba_1, \ba_2, \cdots , \ba_k\}$ of
non-zero vectors of $\nn^m$ (repetition allowed) is called a system of generators on $[m]$ and $\II_\PP$ is the ideal of $\kk[\xx]$ generated by $\xx^{\ba_1},\xx^{\ba_2},\cdots,\xx^{\ba_k}$, i.e., $
\II_\PP=\langle\xx^{\ba_1}, \xx^{\ba_2}, \cdots, \xx^{\ba_k}\rangle$.
Two systems of generators $\PP$ and $\QQ$  are equivalent if they generate the same ideal of $\kk[\xx]$,
i.e., $\II_\PP=\II_\QQ$.}\vspace{2mm}

  {\bf Definition 2.2} {\it For a system of generators $\PP=\{\ba_1, \ba_2, \cdots , \ba_k\}$, the Taylor
DGA algebra $(T(\PP),d)$ over $\kk[\xx]$ is defined as follows. Let $\Lambda(\PP)$ be the exterior algebra over $\kk$ generated by $\PP$, i.e., it is an algebra with unit $1$ generated by $\ba_1,\cdots,\ba_k$ with zero relations $\ba_i\ba_i=0$ and $\ba_i\ba_j=-\ba_j\ba_i$ for all $i,j$. Then algebra $\kk[\xx]{\otimes}\Lambda(\PP)$ ($\otimes$ means the tensor product of algebras over $\kk$) has a DGA algebra over $\kk[\xx]$ structure with differential $d$ defined by $d\ba_i=\xx^{\ba_i}$ for all $1\leqslant i\leqslant k$ and the homological degree defined by $|\ba_i|=1$ for all $1\leqslant i\leqslant m$ and $|\xx^\ba|=0$ for all $\xx^\ba\in\kk[\xx]$. Then $(T(\PP),d)$ is a DAG algebra over $\kk[\xx]$ such that $(\kk[\xx]{\otimes}\Lambda(\PP),d)$ is a sub-DGA of $(T(\PP),d)$ and $T(\PP)$ is generated by $\ba_{i_1,\cdots,i_s}$ ($1\leqslant i_j\leqslant k$) with zero relations $\ba_{i_1}\cdots\ba_{i_s}=\xx^{\ba} \ba_{i_1,\cdots,i_s}$, where $\ba=\ba_{i_1}{+}\cdots{+}\ba_{i_s}{-}[\ba_{i_1},\cdots,\ba_{i_s}]$.
}\vspace{2mm}

By definition, the generators $a_{i_1,\cdots,i_s}$ of $T(\PP)$ satisfies the following properties.

1) $\ba_{i_1,\cdots,i_s}=0$ if $i_u=i_v$ for some $u\neq v$.

2) $\ba_{i_1,\cdots,i_u,\cdots,i_v,\cdots,i_s}=-\ba_{i_1,\cdots,i_v,\cdots,i_u,\cdots,i_s}$.

3) $\ba_{i_1,\cdots,i_s}\ba_{j_1,\cdots,j_t}=\xx^{[\ba_{i_1},\cdots,\ba_{i_s}]+[\ba_{j_1},\cdots,\ba_{j_t}]
-[\ba_{i_1},\cdots,\ba_{i_s},\ba_{j_1},\cdots,\ba_{j_t}]} \ba_{i_1,\cdots,i_s,j_1,\cdots,j_t}$.\vspace{2mm}

{\bf Theorem 2.3} {\it For a system of generators $\PP=\{\ba_1, \ba_2, \cdots , \ba_k\}$,\vspace{1.2mm} \\
\hspace*{40mm}$H_0(T(\PP))=\kk[\xx]/\II_{\PP}$, $H_s(T(\PP))=0$ if $s>0$.\vspace{1.5mm}\\
Thus, $(T(\PP),d)$ is a free resolution of the ring $\kk[\xx]/\II_{\PP}$ as an algebra over $\kk[\xx]$, i.e., the product map of $T(\PP)$ is naturally a chain complex homomorphism $(T(\PP){\otimes}T(\PP),d{\otimes}d)\to (T(\PP),d)$ such that the induced homology homomorphism is just the product map of the ring $\kk[\xx]/\II_{\PP}$.}

\begin{proof} We use induction on $k$, the number of monomials of $\PP$ to prove that $H_0(T(\PP))=\kk[\xx]/\II_{\PP}$ and $H_s(T(\PP))=0$ if $s\neq 0$.

By Leibniz formula, $d(\ba_{i_1}\cdots\ba_{i_s})=\Sigma_{u=1}^s(-1)^{u-1}\xx^{\ba_{i_u}}\ba_{i_1}\cdots\wh \ba_{i_{u}}\cdots\ba_{i_s}$, so
\[d\ba_{i_1,\cdots,i_s}=\Sigma_{u=1}^s(-1)^{u-1}
\xx^{[\ba_{i_1},\cdots,\ba_{i_s}]-[\ba_{i_1},\cdots,\wh \ba_{i_u},\cdots,\ba_{i_s}]}\ba_{i_1,\cdots,\wh i_{u},\cdots,i_s}.\]
Forget the DGA structure of $(T(\PP),d)$ and the above formula makes $(T(\PP),d)$ a chain complex of free modules over $\kk[\xx]$ generated by $1$ and $\ba_{i_1,\cdots,i_s}$ for all $1\leqslant i_1<\cdots<i_s\leqslant k$. If the system of generators has only $1$ monomial, the conclusion is obvious. Suppose the conclusion holds for all system of generators with less than $k$ ($k{>}1$) monomials. Then for $\PP=\{\ba_1, \cdots, \ba_{k}\}$, define $\PP_1=\{\ba_1,\cdots,\ba_{k-1}\}$. From the differential formula we obviously have that $(T(\PP_1),d)$ is a chain subcomplex of $(T(\PP),d)$ and the quotient chain complex is isomorphic to $(\Sigma T(\PP_2),d)$ with $\PP_2{=}\{\ba'_1=[\ba_1,\ba_k],\cdots,\ba'_{k-1}{=}[\ba_{k-1},\ba_k]\}$, where $\Sigma$ means uplifting the degree by $1$ and the isomorphism is defined by corresponding $\ba_{i_1,\cdots,i_s,k}$ in $T(\PP)/T(\PP_1)$ to $\ba'_{i_1,\cdots,i_s}$ in $\Sigma T(\PP_2)$. By the induction hypothesis, $H_0(T(\PP_1))=\kk[\xx]/\II_{\PP_1}$, $H_0(T(\PP_2))=\kk[\xx]/\II_{\PP_2}$ and $H_s(T(\PP_i))=0$ if $s>0$. So from the long exact sequence of homology groups
\[\cdots\to H_1(T(\PP_2))\to H_1(T(\PP_1))\to H_1(T(\PP))\to H_0(T(\PP_2))\to H_0(T(\PP_1))\to H_0(T(\PP))\to 0\]
we have that $H_0(T(\PP))=\kk[\xx]/\II_{\PP}$ and $H_s(T(\PP))=0$ if $s\neq 0$. Thus, $(T(\PP),d)$ is a free resolution of the ring $\kk[\xx]/\II_{\PP}$. By K$\ddot{\rm u}$nneth Theorem, $(T(\PP){\otimes}T(\PP),d{\otimes}d)$ is a free resolution of the ring $\kk[\xx]/\II_{\PP}{\otimes}\kk[\xx]/\II_{\PP}$. The DGA structure of $(T(\PP),d)$ implies that the correspondence $a{\otimes}b\to ab$ is a chain complex homomorphism  $(T(\PP){\otimes}T(\PP),d{\otimes}d)\to (T(\PP),d)$. So we need only prove that $(T(\PP),d)$ is a well-defined DGA algebra over $\kk[\xx]$. This is in fact natural by the DGA algebra structure of $(\kk[\xx]{\otimes}\Lambda(\PP),d)$. Precisely, from the formula
\[d\big((\ba_{i_1}\cdots \ba_{i_s})(\ba_{j_1}\cdots \ba_{j_t})\big)=\big(d(\ba_{i_1}\cdots \ba_{i_s})\big)
\big(\ba_{j_1}\cdots \ba_{j_t}\big){+}(-1)^s\big(\ba_{i_1}\cdots \ba_{i_s}\big)
\big(d(\ba_{j_1}\cdots \ba_{j_t})\big)\]
we have that
\begin{eqnarray*}&&d(\ba_{i_1,\cdots,i_s}\ba_{j_1,\cdots,j_t})\\
&=&\xx^{[\ba_{i_1},\cdots,\ba_{i_s}]+[\ba_{j_1},\cdots,\ba_{j_t}]-[\ba_{i_1},\cdots,\ba_{i_s},\ba_{j_1},\cdots,\ba_{j_t}]}
d(\ba_{i_1,\cdots,i_s,j_1,\cdots,j_t})\\
&=&\sum_{u=1}^s(-1)^{u-1}\xx^{[\ba_{i_1},\cdots,\ba_{i_s}]+[\ba_{j_1},\cdots,\ba_{j_t}]
-[\ba_{i_1},\cdots,\wh\ba_{i_u},\cdots,\ba_{i_s},\ba_{j_1},\cdots,\ba_{j_t}]}
\ba_{i_1,\cdots,\wh i_{u},\cdots,i_s,j_1,\cdots,j_t}\\
&&\!\!\!\!\!\!+\sum_{u=1}^t(-1)^{s+u-1}\xx^{[\ba_{i_1},\cdots,\ba_{i_s}[+[\ba_{j_1},\cdots,\ba_{j_t}]
-[\ba_{i_1},\cdots,\ba_{i_s},\ba_{j_1},\cdots,\wh\ba_{i_u},\cdots,\ba_{j_t}]}
\ba_{i_1,\cdots,i_s,j_1,\cdots,\wh i_{u},\cdots,j_t}\\
&=&\sum_{u=1}^s(-1)^{u-1}\xx^{[\ba_{i_1},\cdots,\ba_{i_s}]-[\ba_{i_1},\cdots,\wh\ba_{i_u},\cdots,\ba_{i_s}]}
\ba_{i_1,\cdots,\wh i_{u},\cdots,i_s}\ba_{j_1,\cdots,j_t}\\
&&\!\!\!\!\!\!+\sum_{u=1}^t(-1)^{s+u-1}\xx^{[\ba_{j_1},\cdots,\ba_{j_t}]-[\ba_{j_1},\cdots,\wh\ba_{i_u},\cdots,\ba_{j_t}]}
\ba_{i_1,\cdots,i_s}\ba_{j_1,\cdots,\wh i_{u},\cdots,j_t}\\
&=&\big(d\ba_{i_1,\cdots,i_s}\big)\big(\ba_{j_1,\cdots,j_t}\big)+(-1)^s \big(\ba_{i_1,\cdots,i_s}\big)\big(d\ba_{j_1,\cdots,j_t}\big)
\end{eqnarray*}
\end{proof}
\vspace{2mm}

For every $\sigma\in 2^{[m]}$, there is a unique vector $(a_1,\cdots,a_i)\in\nn^m$ defined by that $a_i=1$ if $i\in\sigma$ and $a_i=0$ if $i\not\in\sigma$. In this way, we still denote the vector $(a_1,\cdots,a_m)$ by $\sigma$ and regard $2^{[m]}$ as a subset of $\nn^m$ consisting of all those vectors with every coordinate either $0$ or $1$. For $\sigma_1,\cdots,\sigma_s\in 2^{[m]}\subset \nn^m$, $[\sigma_1,\cdots,\sigma_s]=\sigma_1{\cup}\cdots{\cup}\sigma_s$.

By Theorem 1.7 in \cite{Mil} , there is a bijection between square free ideals and simplicial complexes. A monomial ideal generated by $\xx^{\ba_1},\cdots,\xx^{\ba_k}$ is square free if and only if every $\ba_i\in 2^{[m]}$. So we have the following definition.\vspace{1mm}

{\bf Definition 2.4} {\it A simplicial complement $\PP=\{\sigma_1,\cdots,\sigma_k\}$ is a system of generators on $[m]$ such that every $\sigma_i\in 2^{[m]}\subset \nn^m$. The unique simplicial complex $K$ with Stanley-Reisner ring $\kk[\xx]/\II_{\PP}$ is called the simplicial complex associated to $\PP$ and $\PP$ is called the simplicial complement of $K$. The Stanley-Reisner face ring $\kk[\xx]/\II_{\PP}$ of $K$ is denoted by $\kk(K)$.}\vspace{2mm}

{\bf Definition 2.5} {\it Let $\PP=\{\sigma_1,\cdots,\sigma_k\}$ be a simplicial complement on $[m]$ of $K$. The DGA algebra $(\o\Lambda(\PP),d)$ over $\kk$ is defined as follows. $\o\Lambda(\PP)$ is generated as an algebra by elements of the form $\sigma_{i_1}\o\wedge\cdots\o\wedge\sigma_{i_s}$ ($\o\wedge$ is just a symbol but not a product!) satisfying the following properties.

1) $\sigma_{i_1}\o\wedge\cdots\o\wedge\sigma_{i_s}=0$ if $i_u=i_v$ for some $u\neq v$.

2) $\cdots\sigma_{i_u}\o\wedge\cdots\o\wedge\sigma_{i_v}\cdots=-\cdots\sigma_{i_v}\o\wedge\cdots\o\wedge\sigma_{i_u}\cdots$, where the omitted part remain unchanged.

3) $(\sigma_{i_1}\o\wedge\cdots\o\wedge\sigma_{i_s})(\sigma_{j_1}\o\wedge\cdots\o\wedge\sigma_{j_t})=
\left\{\begin{array}{cl}
\sigma_{i_1}\o\wedge\cdots\o\wedge\sigma_{i_s}\o\wedge\sigma_{j_1}\o\wedge\cdots\o\wedge\sigma_{j_t} &
{\rm if}\, \sigma_{i_u}{\cap}\sigma_{j_v}=\phi\,{\rm for\,all}\,u,v;\\
0 & otherwise.
       \end{array}
\right.$

The differential $d$ satisfies
\[d(\sigma_{i_1}\o\wedge\cdots\o\wedge\sigma_{i_s})={\sum}_u(-1)^{u}\sigma_{i_1}\o\wedge\cdots\o\wedge\wh \sigma_{i_u}\o\wedge\cdots\o\wedge\sigma_{i_s},\]
where the sum is taken over all $u$ such that $\sigma_{i_1}{\cup}\cdots{\cup}\wh \sigma_{i_u}{\cup}\cdots{\cup}\sigma_{i_s}=\sigma_{i_1}{\cup}\cdots{\cup}\sigma_{i_s}$.

For $\sigma\in K^c$, let $\o\Lambda^{s,\sigma}(\PP)$ be the submodule of $\o\Lambda^s(\PP)$ generated by all $\sigma_{i_1}\o\wedge\cdots\o\wedge\sigma_{i_s}$ such that $\sigma=\sigma_{i_1}{\cup}\cdots{\cup}\sigma_{i_s}$. It is obvious that $(\o\Lambda^{*,\sigma}(\PP),d)$ is a chain subcomplex of $(\o\Lambda(\PP),d)$ such that\\
\hspace*{36mm}$\o\Lambda^{*,\sigma}(\PP)\o\Lambda^{*,\sigma'}(\PP)\subset
\left\{\begin{array}{cl}
\o\Lambda^{*,\sigma\cup\sigma'}(\PP)&
{\rm if}\, \sigma{\cap}\sigma'=\phi;\\
0 & otherwise.
       \end{array}
\right.$}\vspace{2mm}

{\bf Theorem 2.6} Let $\PP$ be a simplicial complement on $[m]$ of $K$. There is a direct sum decomposition
\[Tor^{\kk[\xx]}_*(\kk(K),\kk)=\kk\oplus\left(\oplus_{\sigma\in K^c}H_*(\o\Lambda^{*,\sigma}(\PP),d)\right),\]
such that $Tor^{\kk[\xx]}_0(\kk(K),\kk)=\kk$ and for all $i>0$ and $\sigma\in K^c$,
\[Tor^{\kk[\xx]}_{*,\sigma}(\kk(K),\kk)=H_i(\o\Lambda^{*,\sigma}(\PP),d)=\w H_{i-2}({\rm link}_{K^*}\o\sigma;\kk)
=\w H^{|\sigma|-i-1}(K{\cap}\sigma;\kk),\]
where $\kk$ denotes the trivial chain complex.

\begin{proof} By Theorem 2.3, for any module $N$ over $\kk[\xx]$, $Tor^{\kk[\xx]}_*(\kk(K),\kk)=(T(\PP){\otimes}_{\kk[\xx]}N,d{\otimes}1_N)$ and if $N$ is an algebra over $\kk[\xx]$, then $(T(\PP){\otimes}_{\kk[\xx]}N,d{\otimes}1_N)$ is naturally a DGA with product defined by that for all $n_1,n_2\in N$,
\[(\ba_{i_1,\cdots,i_s}{\otimes}n_1)(\ba_{j_1,\cdots,j_t}{\otimes}n_2)=\ba_{i_1,\cdots,i_s,j_1,\cdots,j_t}
{\otimes}\xx^\ba n_1n_2,\]
where $\ba=[\ba_{i_1},\cdots,\ba_{i_s}]{+}[\ba_{j_1},\cdots,\ba_{j_t}]{-}
[\ba_{i_1},\cdots,\ba_{i_s},\ba_{j_1},\cdots,\ba_{j_t}]$. Take $N=\kk$, the trivial algebra over $\kk[\xx]$ and we have
$(T(\PP){\otimes}_{\kk[\xx]}\kk,d)=(\o\Lambda(\PP),d)$, with $\ba_{i_1,\cdots,i_s}{\otimes}1$ denoted by $\sigma_{i_1}\o\wedge\cdots\o\wedge\sigma_{i_s}$. So $Tor^{\kk[\xx]}_*(\kk(K),\kk)=H_*(\o\Lambda(\PP),d)$ and the DGA structure of $(\o\Lambda(\PP),d)$ is induced from the product of $\kk$. Thus, the direct sum decomposition exists and the isomorphism is an algebra isomorphism. The direct proof that $H_i(\o\Lambda^{*,\sigma}(\PP),d)=\w H_{i-2}({\rm link}_{K^*}\o\sigma;\kk)$ is too long and we will give it in a proceeding paper.
\end{proof}
\vspace{2mm}

{\bf Example 2.7} Compute the algebra $Tor^{\kk[\xx]}_{*}(\kk(K), \kk)$, where $K$ is as shown in Figure 1.
\begin{center}
\setlength{\unitlength}{0.3mm}
\begin{picture}(210,120)(0,-10)
\put(0,0){\line(1,0){200}}
\put(0,0){\line(1,1){100}}
\put(200,0){\line(-1,1){100}}
\put(66.6,33.3){\line(1,0){66.6}}
\put(66.6,33.3){\line(1,2){33.3}}
\put(100,100){\line(1,-2){33.3}}
\put(0,0){\line(2,1){66.6}}
\put(200,0){\line(-2,1){66.6}}
\put(-7,-5){\mbox{1}}
\put(205,-5){\mbox{2}}
\put(66.6,20){\mbox{4}}
\put(130,20){\mbox{5}}
\put(98,103){\mbox{3}}
\put(80,-15){\mbox{Figure 1}}
\put(10,10){\line(0,-1){5}}
\put(20,20){\line(0,-1){10}}
\put(30,30){\line(0,-1){15}}
\put(40,40){\line(0,-1){20}}
\put(50,50){\line(0,-1){25}}
\put(60,60){\line(0,-1){30}}
\put(70,70){\line(0,-1){30}}
\put(80,80){\line(0,-1){20}}
\put(90,90){\line(0,-1){9}}
\put(110,90){\line(0,-1){9}}
\put(120,80){\line(0,-1){20}}
\put(130,70){\line(0,-1){30}}
\put(140,60){\line(0,-1){30}}
\put(150,50){\line(0,-1){25}}
\put(160,40){\line(0,-1){20}}
\put(170,30){\line(0,-1){15}}
\put(180,20){\line(0,-1){10}}
\put(190,10){\line(0,-1){5}}
\end{picture}
\end{center}\vspace{5mm}

$K^c=\{\{1,5\},\{2,4\},\{1,2,3\},\{1,2,4\},\{1,2,5\},\{1,3,5\},\{1,4,5\},\{2,3,4\},\{2,4,5\},\{3,4,5\}\}$. So we may take $\PP=\left\{\sigma_1=\{1,5\},\ \sigma_2=\{2,4\},
\ \sigma_3=\{1,2,3\},\ \sigma_4=\{3,4,5\}\right\}$ to be the simplicial complement of $K$.

We directly compute the chain complex $(\o\Lambda(\PP),d)$. By definition, all the non-trivial differentials are
\begin{align*}
d(\sigma_1{\o\wedge}\sigma_2{\o\wedge}\sigma_3{\o\wedge}\sigma_4)= & -\sigma_2{\o\wedge}\sigma_3{\o\wedge}\sigma_4 +\sigma_1{\o\wedge}\sigma_3{\o\wedge}\sigma_4
-\sigma_1{\o\wedge}\sigma_2{\o\wedge}\sigma_4 +\sigma_1{\o\wedge}\sigma_2{\o\wedge}\sigma_3;\\
d(\sigma_1{\o\wedge}\sigma_3{\o\wedge}\sigma_4) = & -\sigma_3{\o\wedge}\sigma_4;\\
d(\sigma_2{\o\wedge}\sigma_3{\o\wedge}\sigma_4) = & -\sigma_3{\o\wedge}\sigma_4.
\end{align*}
Thus the $Tor^{\kk[\xx]}_{*}(\kk(K), \kk)$ is a vector space over $\kk$ with a base
\[
\left\{\begin{array}{ll}
         \sigma_1{\o\wedge}\sigma_2{\o\wedge}\sigma_4\\
         \sigma_1{\o\wedge}\sigma_2{\o\wedge}\sigma_3
         \end{array}
        \begin{array}{ll}
         \sigma_1{\o\wedge}\sigma_2\\
         \sigma_1{\o\wedge}\sigma_3\\
         \sigma_1{\o\wedge}\sigma_4\\
         \sigma_2{\o\wedge}\sigma_3\\
         \sigma_2{\o\wedge}\sigma_4
        \end{array}
        \begin{array}{ll}
         \sigma_1\\
         \sigma_2\\
         \sigma_3\\
         \sigma_4
        \end{array}
        \begin{array}{ll}
         1
        \end{array}
\right\}.
\]

Thus, $Tor_{*}^{\kk[\xx]}(\kk(K), \kk)$ is a free $\kk$-module with Poincar\`{e}
series $1+4x+5x^2+2x^3$. The algebraic structure of $Tor_{*}^{\kk[\xx]}(\kk(K), \kk)$ is very easy by the product of $\Lambda(\PP)$. The only non-trivial product is $\sigma_1\sigma_2=\sigma_1{\o\wedge}\sigma_2$.\vspace{2mm}

This example shows the advantage of the chain complex $(\o\Lambda(\PP),d)$. Usually, the computation of $Tor_*(\kk(K),\kk)$ is by Hochster formula, since the  Koszul complex tensor with $\kk(K)$ is still a module over $\kk(K)$ and the direct computation the complex is formidable. But in this example, we may directly compute the DGA algebra $(\o\Lambda(\PP),d)$ without using Hochster formula.\vspace{2mm}

{\bf Definition 2.8} {\it Let $\PP=\{\sigma_1, \sigma_2, \cdots, \sigma_s\}$ be a simplicial complement on $[m]$ of $K$.
For $\omega\subset[m]$, the $\omega$-compression of $\PP$  is defined by
\[
E_\omega\PP=\{\sigma_1{\setminus}\omega, \sigma_2{\setminus}\omega, \cdots, \sigma_s{\setminus}\omega\}.
\]}

    For an arbitrary simplex $\omega\in K$, define its {\it link} and {\it star} to be the
subcomplexes
\begin{align*}
  & \mbox{star}_K\omega = \{\tau\in K|\,\omega{\cup}\tau\in K\},
  & \mbox{link}_K\omega = \{\tau\in K|\,\omega{\cup}\tau\in K, \omega{\cap}\tau=\phi\}={\rm star}_K\omega\cap([m]{\setminus}\omega).
\end{align*}

{\bf Theorem 2.9} {\it Let $\PP$ be a simplicial complement of the simplicial complex $K$. Then for any $\omega\subset[m]$, $E_\omega\PP$ is the simplicial complement of the simplicial complex ${\rm star}_{K}\omega$. Thus
for any $\sigma\in {\rm star}_K\omega$ and $i>0$,\[
H_{i,\sigma}(\o\Lambda^{*,\sigma}(E_\omega\PP), d)=
  Tor^{\kk[\xx]}_{i,\sigma}(\kk({\rm star}_{K}\omega), \kk),
\]
and
\[
H_{q,[m]\setminus \omega}(\o\Lambda^{*,[m]\setminus \omega}(E_\omega\PP), d)=
 \w{H}^{m-|\omega|-q-1}({\rm link}_{K}\omega, \kk).
\]}

\begin{proof}
We may suppose $\PP=\{K^c\}$, the set of non-faces of $K$ (in any order). Then
\[({\rm star}_K\omega)^c=\{\tau\,|\,\omega{\cup}\tau\not\in K\}=\{\sigma{\setminus}\omega\,|\,\sigma\in K^c\}
=E_{\omega}\PP.\]
So by Theorem 2.6 and Hochster formula $H_{q,\o\omega}(\o\Lambda^{*,\o\omega}(E_\omega\PP), d)=
 \w{H}^{m-|\omega|-q-1}({\rm star}_K\omega{\cap}\o\omega, \kk)= \w{H}^{m-|\omega|-q-1}({\rm link}_{K}\omega, \kk)$, where $\o\omega=[m]{\setminus}\omega$.
\end{proof}

\section{The cohomology of generalized moment-angle complexes}

    In this section we consider the cohomology module of the generalized moment-angle complex.
Recall from \cite{BB} \cite{BBCG08},\vspace{2mm}

  {\bf Definition \nr} {\it Let $({X}, {A})=\{(X_i, A_i, x_i)|i\in[m]\}$ be
a set of pointed $CW$-pairs $(X_i,A_i,x_i)$ ($(A_i,x_i)$ is a pointed CW-subcomplex of $(X_i,x_i)$) and $K$ be an abstract simplicial complex. The generalized moment-angle complex determined by $({X}, {A})$ and
$K$ denoted by $\ZZ_{K}({X}, {A})$ is defined as follows. For every $\omega\in K$, let
$D(\omega)= Y_1{\times} Y_2 {\times} \cdots {\times} Y_m $, where $Y_i=X_i$ if $i\in\omega$ and $Y_i=A_i$ if $i\not\in\omega$. Then the generalized moment-angle complex is
\[
\ZZ_{K}({X}, {A})=\bigcup_{\omega\in K} D(\omega).
\]

    Let $Y_1{\wedge }Y_2{ \wedge} \cdots {\wedge} Y_m $ be the smash product given by
the quotient space
\[
Y_1{\times} Y_2{\times}\cdots{\times} Y_m/S(Y_1{\times} Y_2{\times}\cdots{\times} Y_m)
\]
where $S(Y_1{\times} Y_2{\times}\cdots{\times} Y_m)$ is the subspace of the product with at least
one coordinate given by the base-point $x_i\in Y_i$. The {\it generalized smash moment-angle
complex} is defined to be the image of $\ZZ_{K}({X}, {A})$
in $X_1{\wedge} X_2{\wedge}\cdots{\wedge} X_m$, i.e.,
\[
\wh{\ZZ}_{K}({X}, {A})=\bigcup_{\omega\in K} \wh{D}(\omega),
\]
where
$\wh{D}(\omega)= Y_1{\wedge}Y_2{\wedge}\cdots{\wedge}Y_m$  with $Y_i=X_i$ if $i\in\omega$ and $Y_i=A_i$ if $i\not\in\omega$.

Given a non-empty subset $I=\{i_1, i_2, \cdots, i_k\}\subset [m]$ and a family of pointed pairs
$({X}, {A})$, define
\[
({X_I}, {A_I})=\{(X_{i_j}, A_{i_j})|i_j\in I\}
\]
which is the subfamily of $({X}, {A})$ determined by $I$.  It is known from
\cite{BB,BBCG08} that
\begin{align}\tag{\nr}
H: \Sigma(\ZZ_{K}({X}, {A}))\longrightarrow &
  \Sigma\left(\bigvee_{I\subset [m]}\wh{\ZZ}_{K\cap I}({X_I}, {A_I})\right)
\end{align}
is a natural pointed homotopy equivalence. Precisely, let
$\omega\subset I$ and
$\wh{D}_I(\omega)=Y_{i_1}{\wedge}Y_{i_2}{\wedge}\cdots{\wedge}Y_{i_k}$ with $Y_{i_j}=X_{i_j}$ if $i_j\in \omega$ and $Y_{i_j}=A_{i_j}$ if $i_j\not\in\omega$. Then
\[
\wh{\ZZ}_{K\cap I}({X_I}, {A_I})=\bigcup_{\omega\in K\cap I}\wh{D}_I(\omega).
\]}

    Associated to a simplicial complex $K$, there is a partial ordered set ({\it poset}) $\overline{K}$
with point $\sigma$ in $\overline{K}$ corresponding to a simplex $\sigma\in K$ and order given by
reverse inclusion of simplices. Thus $\sigma_1\leqslant \sigma_2$ in $\overline{K}$ if and only if
$\sigma_2\subseteq\sigma_1$ in $K$.  Given an $\omega\in\overline{K}$ there are further posets given
by
\begin{align*}
\overline{K}_{<\omega} = & \{\tau\in\overline{K}|\tau<\omega\}=\{\tau\in\overline{K}|\omega\varsubsetneq\tau\}
 & \mbox{and} \\
\overline{K}_{\leqslant\omega} = & \{\tau\in\overline{K}|\tau\leqslant\omega\}
 =\{\tau\in\overline{K}|\omega\subset\tau\}.
\end{align*}

    Given a poset $P$, the {\it order complex} $\Delta(P)$ is the simplicial complex with vertices given by set of
points of $P$ and $k$-simplices given by ordered $(k+1)$-tuples $(p_0, p_1, \cdots, p_k)$ with
$p_0<p_1<\cdots<p_k$.  It follows that $\Delta(\overline{K})_{<\phi}=K'_{\PP}$ is the barycentric subdivision
of $K$.

    Given a simplicial complex $K$, we use $|K|$ to denote its geometric realization.
The symbol $|K|\ast Y$ denotes the join of $K$ and $Y$. If $Y$ is a pointed $CW$-complex, then
$|K|\ast Y$ has the homotopy type of $\Sigma |K|\wedge Y$.\vspace{1mm}

  {\bf Theorem (Bahri, Bendersky, Cohen and Gitler \cite{BB,BBCG08} 2.12)} {\it Let $K$ be an abstract simplicial
complex and $\overline{K}$
its associated poset. Let $({X}, {A})=\{(X_i, A_i, x_i)|i\in [m]\}$ denote m choices of
connected, pointed pairs of $CW$-complexes, with the inclusion $A_i\hookrightarrow X_i$ homotopic to constant for
all $i$. Then there is a homotopy equivalence
\[
\wh{\ZZ}_{K}({X}, {A})\longrightarrow \bigvee_{\omega\in K}
|\Delta(\overline{K}_{<\omega})|\ast\wh{D}(\omega).
\]
}

    From (3.2), we see that if $A_i\hookrightarrow X_i$ are null-homotopic for all $i$, then
\begin{align*}
\Sigma(\ZZ_{K}({X}, {A}) \simeq &
  \Sigma\left(\bigvee_{I\subset [m]} \wh{\ZZ}_{K \cap I}({X_I}, {A_I})\right)\\
   \simeq &  \Sigma\left(\bigvee_{I\subset [m]}\left(\bigvee_{\omega\in K\cap I}
  |\Delta(\overline{K\cap I}_{<\omega})|\ast\wh{D}_I(\omega)\right)\right).\notag
\end{align*}
Fix an $\omega\in K$, it is easy to see that $\omega\in K\cap I$ if and only if
$\omega\subset I$. Thus $\Sigma(\ZZ_{K}({X}, {A})$ has the homotopy type of
\begin{align}\notag
\left(\bigvee_{I\subset [m]}\left(\bigvee_{\omega\in K\cap I}
  |\Delta(\overline{K\cap I}_{<\omega})|\ast\wh{D}_I(\omega)\right)\right) = &
\left(\bigvee_{\omega\in K}\left(\bigvee_{\omega \subset I\subset [m]}
  |\Delta(\overline{K\cap I}_{<\omega})|\ast\wh{D}_I(\omega)\right)\right)\\
 \simeq & \left(\bigvee_{\omega\in K}\left(\bigvee_{\omega \subset I\subset [m]}
  \Sigma|\Delta(\overline{K\cap I}_{<\omega})|\wedge\wh{D}_I(\omega)\right)\right)\tag{\nr}
\end{align}

  {\bf Remark \nr} If $\omega$ is a maximum face of $K\cap I$
in the sense that $\omega\in K\cap I$ but it is not a proper subset of any other simplices,
then $\overline{K\cap I}_{<\omega}=\{\tau\in \overline{K\cap I}|\omega\varsubsetneq\tau\}=\varnothing$.
The simplicial complex $\Delta(\overline{K\cap I}_{<\omega})=\{\phi\}$ is the empty simplicial complex.
Here we use the agreement
\begin{align*}
|\{\phi\}|\ast \wh{D}_I(\omega)= & \wh{D}_I(\omega) & \mbox{and} & &
\Sigma|\{\phi\}|\wedge \wh{D}_I(\omega) = & \wh{D}_I(\omega).
\end{align*}
Combine with the agreement in Theorem 2.8, $\w{H}^{-1}(|\{\phi\}|, \kk)=\kk$, we have
\begin{align*}
\w{H}^0(\Sigma|\{\phi\}|, \kk)= & \kk & \mbox{and} & & \w{H}^*(\Sigma|\{\phi\}|\wedge\wh{D}_I(\omega), \kk)= &
\w{H}^*(\wh{D}_I(\omega), \kk)
\end{align*}

   Consider the reduced cohomology of $\ZZ_{K}({X}, {A})$, we see from (5.3) that:
\begin{align}\tag{\nr}
\w{H}^*(\ZZ_{K}({X}, {A}), \kk)= &
  \bigoplus_{\omega\in K}\left(\bigoplus_{\omega\subset I\subset[m]}
  \w{H}^*\left(\Sigma|\Delta(\overline{K\cap I}_{<\omega})|\wedge\wh{D}_I(\omega), \kk\right)\right).
\end{align}

   Furthermore suppose that there is no
$Tor$ problem in the K\"{u}nneth formulae for the cohomology of
$\Sigma|\Delta(\overline{K\cap I}_{<\omega})|\wedge\wh{D}_I(\omega)$
(For example, take $\kk$ to be a field or suppose that $\w{H}^*(X_i, \kk)$ and $\w{H}^*(A_i, \kk)$ are free
$\kk$-modules for all $i$).
Then from
\begin{align*}
\wh{D}_I(\omega)= & Y_{i_1}\wedge Y_{i_2}\wedge\cdots\wedge Y_{i_k}\simeq
 \left(\bigwedge_{i\in\omega}X_i\right)\bigwedge\left(\bigwedge_{j\in I\setminus\omega} A_j\right),\\
H^*(\wh{D}_I(\omega), \kk) = & \left(\bigotimes_{i\in\omega} \w{H}^*(X_i, \kk)\right) \otimes\left(\bigotimes_{j\in I\setminus\omega}
  \w{H}^*(A_j, \kk)\right)
\end{align*}
which is a free $\kk$-module, we have
\begin{align}\tag{\nr}
 & \w{H}^*\left(\Sigma|\Delta(\overline{K\cap I}_{<\omega})|\wedge\wh{D}_I(\omega), \kk\right)
=  \w{H}^*\left(\Sigma|\Delta(\overline{K\cap I}_{<\omega})|, \kk\right)\otimes\w{H}^*\left(\wh{D}_I(\omega), \kk\right)\\
= & \w{H}^*(\Sigma|\Delta(\overline{K\cap I}_{<\omega})|,\kk)\otimes
  \left(\bigotimes_{i\in\omega} \w{H}^*(X_i, \kk)\right) \otimes\left(\bigotimes_{j\in I\setminus\omega}
  \w{H}^*(A_j, \kk)\right).\notag
\end{align}

 {\bf Theorem \nr} {\it Let $K$ be a (abstract) simplicial complex with simplicial complement $\PP$
and $({X}, {A})=\{(X_i, A_i, x_i)|i\in [m]\}$  a set of pointed $CW$-complex pairs such that the inclusion $A_i\hookrightarrow X_i$ homotopic to constant for
all $i$. Then the cohomology  of $\ZZ_{K}({X}, {A})$ over $\kk$ is isomorphic to
\[
H^*(\ZZ_{K}({X}, {A}), \kk)=
\bigoplus_{\omega\in K}\left(\bigoplus_{\sigma}
  H_{*,\sigma}(\Lambda^{*,\sigma}(E_\omega\PP), d) \otimes
  \left(\bigotimes_{i\in\omega}\w{H}^*(X_i, \kk)\right)\otimes
  \left(\bigotimes_{j\in\sigma}\w{H}^*(A_j, \kk)\right)\right)
\]
as vector spaces over $\kk$,  where
\[
H_{*,\sigma}(\Lambda^{*,\sigma}(E_\omega\PP),d)\cong Tor_{*,\sigma}^{\kk[\xx]}(\kk(\mbox{star}_{K}\omega), \kk)
\]
subject to $\sigma\subset [m]\setminus \omega$.
}

\begin{proof}
    Given an $\omega\subset I$, from its definition we see that the posets
\begin{align*}
\overline{K\cap I}_{\leqslant \omega}= & \{\tau\in K\cap I|\omega\subset \tau\}
  = \{\tau\in K|\omega\subset\tau\subset I\} \hspace{10mm}\mbox{and}\\
\overline{\mbox{link}_{K\cap I}\omega} = & \{\tau'\in K\cap I|\omega\cup\tau'\in K\cap I, \tau'\cap\omega=\phi\}\\
  = & \{\tau'\in K|\omega\cup\tau'\in K, \tau'\cap\omega=\phi,\ \mbox{and}\ \tau'\subset I\}
     =\overline{(\mbox{link}_{K}\omega)\cap I}_{\leqslant \phi}.
\end{align*}
There is a one to one correspondence between the posets $\Psi: \overline{K\cap I}_{<\omega}\longrightarrow
\overline{(\mbox{link}_{K}\omega)\cap I}_{<\phi}$ given by $\Psi(\tau)=\tau\setminus\omega$.
Thus the order complex $\Delta(\overline{K\cap I}_{<\omega})=((\mbox{link}_{K}\omega)\cap I)'$
is the barycentric subdivision of
$(\mbox{link}_{K}\omega)\cap I$ and then
\begin{align}\tag{\nr}
\w{H}^*(\Sigma|\Delta(\overline{K\cap I}_{<\omega})|, \kk)=
\w{H}^{*-1}(|\Delta(\overline{K\cap I}_{<\omega})|, \kk)=
  \w{H}^{*-1}((\mbox{link}_{K}\omega)\cap I, \kk).
\end{align}

Since ${\rm link}_{K}\omega={\rm star}_K\omega\cap([m]{\setminus}\omega)$, we have
\[({\rm link}_{K}\omega)\cap I =
{\rm star}_K\omega\cap([m]{\setminus}\omega)\cap I={\rm star}_K\omega\cap (I{\setminus}\omega).
\]
from Theorem 2.9 and (3.8) we have
\begin{align*}
H_{q,I\setminus\omega}(\Lambda^{*,I\setminus\omega}(E_\omega\PP), d)
 =  & \w{H}^{|I\setminus\omega|-q-1}((\mbox{link}_{K}\omega)\cap I, \kk)
 =  \w{H}^{|I\setminus\omega|-q}(\Sigma|\Delta(\overline{K\cap I}_{<\omega})|, \kk)
\end{align*}

From the definition of $E_\omega\PP=\{\sigma_1{\setminus}\omega, \sigma_2{\setminus}\omega, \cdots,
\sigma_s{\setminus}\omega\}$, we see that  the homology of the simplicial complement
$H_{q,\sigma}(\Lambda^{*,\sigma}(E_\omega\PP),d)$ is concentrated in $\sigma\subset[m]{\setminus}\omega$.
For $\omega\subset I\subset[m]$, denote $\sigma=I{\setminus}\omega$. Then
we have
\[
H_{q,\sigma}(\Lambda^{*,\sigma}(E_\omega\PP), d)=\w{H}^{|\sigma|-q}(\Sigma|\Delta
 (\overline{K\cap (\sigma\cup\omega)}_{<\omega})|, \kk).
\]

  Apply this formula to (3.6) and (3.5) we get
\begin{align*}
 & \w{H}^*(\ZZ_{K}({X}, {A}), \kk)=
  \bigoplus_{\omega\in K}\left(\bigoplus_{\sigma\subset[m]\setminus\omega}
  \w{H}^*\left(\Sigma|\Delta(\overline{K\cap (\omega\cup\sigma)}_{<\omega})|\wedge\wh{D}_{\omega\cup\sigma}(\omega), \kk\right)\right) \\
 = &  \bigoplus_{\omega\in K}\left(\bigoplus_{\sigma\subset[m]\setminus\omega}
 H_{q, \sigma}(\Lambda^{*, \sigma}(E_\omega\PP), d)\otimes
  \left(\bigotimes_{i\in\omega} \w{H}^*(X_i, \kk)\right) \otimes\left(\bigotimes_{j\in \sigma}
  \w{H}^*(A_j, \kk)\right)\right)
\end{align*}
where $I=\omega\cup\sigma\not=\phi$.

    The theorem follows from $H^0(\ZZ_{K}({X}, {A}), \kk)=\kk$ and the agreement by taking
$I=\omega=\sigma=\phi$ to be the empty set
\[
H_{0,\phi}(\Lambda^{*,\phi}[E_\phi\PP], d)
  \otimes\left(\bigotimes_{i\in\phi}\w{H}^*(X_i, \kk)\right)
  \otimes\left(\bigotimes_{j\in\phi}\w{H}^*(A_j, \kk)\right)=\kk.
\]
\end{proof}

   {\bf Remark} For $\omega\not\in K$, the empty set $\phi$ is an element of $E_\omega\PP$ and
$H_{*,*}(\o\Lambda^{*,*}(E_\omega\PP), d)\equiv 0$. Thus Theorem 3.7 could be written as
\begin{align}\notag
  & H^*(\ZZ_{K}({X}, {A}), \kk)\\
= & \bigoplus_{\omega\in[m]}\left(\bigoplus_{\sigma\in[m]}
  H_{*,\sigma}(\o\Lambda^{*,\sigma}(E_\omega\PP), d) \otimes
  \left(\bigotimes_{i\in\omega}\w{H}^*(X_i, \kk)\right)\otimes
  \left(\bigotimes_{j\in\sigma}\w{H}^*(A_j, \kk)\right)\right)\tag{\nr}
\end{align}

  {\bf Corollary \nr} {\it If all the $A_i$ are contractible, then
\[
H^*(\ZZ_{K}({X}, {A}), \kk)=
 \bigoplus_{\omega\in K}\left(\bigotimes_{i\in \omega} \w{H}^*(X_i, \kk)\right)
\]
}

\begin{proof}
If all the $A_i$ are contractible then $\displaystyle{\bigotimes_{j\in\tau}\w{H}^*(A_j, \kk)=0}$ for any non-empty set $\tau$.
The corollary follows from $H_{0,\phi}(\o\Lambda^{*,\phi}(E_\omega\PP),d)=\kk$ if $\omega\in K$.
\end{proof}

 {\bf Corollary \nr} {\it If all the $X_i$ are contractible, then
\begin{align*}
H^*(\ZZ_{K}({X}, {A}), \kk)= &
 \bigoplus_{\sigma}\left(H_{*,\sigma}(\o\Lambda^{*,\sigma}(\PP), d)
 \otimes\left(\bigotimes_{j\in\sigma}\w{H}^*(A_j, \kk)\right)\right).
\end{align*}
Furthermore take $X_i=D^2$ and $A_i=S^1$ for all $i$,
\begin{align*}
H^{r}(\ZZ_{K}, \kk)=H^{2|\sigma|-q}(\ZZ_{K}({D^2}, {S^1}), \kk)= &
 \bigoplus_{2|\sigma|-q=r}H_{q,\sigma}(\o\Lambda^{*,\sigma}(\PP), d),
\end{align*}
where $H_{q,\sigma}(\o\Lambda^{*,\sigma}(\PP), d)=Tor_{q,\sigma}^{\kk[\xx]}(\kk(K), \kk)$.
}

\begin{proof}
    If all the $X_i$ are contractible, then $\displaystyle{\bigotimes_{i\in\omega}\w{H}^*(X_i, \kk)=0}$ for any
non-empty set $\omega$. The corollary follows from
\[
H^*(\ZZ_{K}({X}, {A}), \kk)=
 \bigoplus_{\sigma}H_{*,\sigma}(\o\Lambda^{*,\sigma}[E_\phi\PP], d)
 \otimes\left(\bigotimes_{i\in\phi}\w{H}^*(X_i,\kk)\right)\otimes
  \left(\bigotimes_{j\in\sigma}\w{H}^*(A_j, \kk)\right)
\]
and the agreement $\displaystyle{\bigotimes_{i\in\phi}\w{H}^*(X_i,\kk)=\kk}$.
\end{proof}

  {\bf Proposition \nr} {\it Let $X_i=S^2$ and $A_i=S^1$ all $i$. Then
\[
H^{r}(\ZZ_{K}({S^2}, {S^1}), \kk) =
  \bigoplus_{\begin{array}{c}\omega\in K\\ 2|\omega|+2|\sigma|-q=r\end{array}}
  \left(\bigoplus_{\sigma}H_{q,\sigma}(\o\Lambda^{*,\sigma}(E_\omega\PP), d)\right).
\]
The isomorphism is an algebra isomorphism by Theorem 7.6 of \cite{PB}.}

\begin{proof}
  Noticed that $\wh{D}_{\omega\cup\sigma}(\omega)$ is the $2|\omega|+|\sigma|$ sphere if $\omega\cap\sigma=\phi$,
$$\wh{D}_{\omega\cup\sigma}(\omega)
 =\left(\bigwedge_{i\in\omega}S^2\right)\wedge\left(\bigwedge_{j\in\sigma}S^1\right)
 =S^{2|\omega|+|\sigma|},
$$
we see from (3.5) that
$$H^*(\ZZ_{K}({X}, {A}), \kk)=
  \bigoplus_{\omega\in K}\left(\bigoplus_{\sigma\subset[m]\setminus\omega}
  \w{H}^*\left(\Sigma^{2|\omega|+|\sigma|+1}|\Delta(\overline{K\cap (\omega\cup\sigma)}_{<\omega}|, \kk\right)\right)
$$
including $\omega=\sigma=\phi$. The result follows from
$$\w{H}^{|\sigma|-q}\left(\Sigma|\Delta(\overline{K\cap (\omega\cup\sigma}_{<\omega})|, \kk\right)
 =H_{q,\sigma}(\o\Lambda^{*,\sigma}(E_\omega\PP), d).$$
\end{proof}\vspace{2mm}

  We finish this paper by giving an example.\vspace{1mm}

{\bf Example \nr} Let $m=6$ and the simplicial complement $\PP=\{\sigma_1=\{1,2\}, \sigma_2=\{3,4\}, \sigma_3=\{5,6\}\}$.
The corresponding simplicial complex is a triangulation of the sphere $S^2$ (see Finger 2)

\begin{center}
\setlength{\unitlength}{0.5mm}
\begin{picture}(100,60)(0,0)
\Thicklines
\drawline(0,0)(100,0)
\drawline(0,0)(0,50)
\drawline(0,0)(25,25)
\drawline(0,50)(25,25)
\drawline(0,50)(100,50)
\drawline(25,25)(125,25)
\drawline(100,0)(125,25)
\drawline(100,50)(125,25)
\thinlines
\drawline(0,0)(125,25)
\drawline(25,25)(100,50)
\dottedline{2}(0,50)(100,0)
\Thicklines
\dottedline{2}(100,0)(100,50)
\put(-6,-1){1}
\put(-6,50){3}
\put(15,22){5}
\put(103,50){2}
\put(127,22){4}
\put(105,-1){6}
\put(35,-20){Finger 2}
\end{picture}
\end{center}
\vspace{15mm}

    The cohomology ring of $\ZZ_K=\ZZ_{K}({D^2}, {S^1})$ can be easily obtained from
the homology of the simplicial complement  $\PP=\{\sigma_1=\{1,2\}, \sigma_2=\{3,4\},$ $\sigma_3=\{5,6\}\}$ and
\[
H_{*,*}(\o\Lambda^{*,\tau}(\PP), d)\cong\o\Lambda^{*,*}(\{\sigma_1, \sigma_2, \sigma_3\}).
\]
Since the algebra $\o\Lambda^{*,*}(\{\sigma_1, \sigma_2, \sigma_3\})$ is just the exterior algebra $\Lambda(\sigma_1,\sigma_2,\sigma_3)$ over $\kk$ generated by $\PP$, so is the cohomology ring $H^*(\ZZ_K, \kk)$.
The Poincar\`{e} series of $H^*(\ZZ_K, \kk)$ is  $1+3x^3+3x^6+x^9$ and the
total Betti number of $\ZZ_K$ is $8$.

Now  consider the cohomology of $\ZZ_{K}({S^2}, {S^1})$. We start from
computing the homology of the simplicial complement $E_\omega\PP$ with $\omega\in K$:
\begin{enumerate}
\item Take $\omega=\phi$,  the homology
of the simplicial complement $E_\phi\PP$ is
$\Lambda(\sigma_1, \sigma_2, \sigma_3).$
Thus the graded vector space
\[
\bigoplus_{\sigma}H_{*,\sigma}(\o\Lambda^{*,\sigma}(E_\phi\PP), d)=\Lambda^{*,*}(\{\sigma_1,\sigma_2,\sigma_3\})
\]
has Poincar\`{e} series $1+3x^3+3x^6+x^9$.
\item Take $\omega=\{1\}$ (similarly, take any of its 6 $0$-simplex $\omega=\{i\}$, $i\in[6]$),
the homology
of the simplicial complement $E_{\{i\}}\PP$ is $\o\Lambda^{*,*}[E_{\{i\}}\PP]$.
Thus the graded vector space
\[
\displaystyle{\bigoplus_{\sigma}H_{*,\sigma}(\o\Lambda^{*,\sigma}(E_{\{i\}}\PP), d)}
\]
has Poincar\`{e} series $x^2(1+x+2x^3+2x^4+x^6+x^7)$.
\item Take $\omega=\{1,3\}$ and similarly any of its 12 $1$-simplex $\omega=\{i_1, i_2\}$, the homology
of the corresponding simplicial complement is $\o\Lambda^{*,*}(E_{\{i_1,i_2\}}\PP)$. The graded vector space
\[
\displaystyle{\bigoplus_{\sigma}H_{*,\sigma}(\o\Lambda^{*,\sigma}(E_{\{i_1, i_2\}}\PP), d)}
\]
is the free $\kk$-module with Poincar\`{e} series $x^4(1+2x+x^2+x^3+2x^4+x^5)$.
\item Take $\omega=\{1,3,5\}$ and similarly any of its 8 $2$-simplex $\omega=\{i_1, i_2,i_3\}$, the homology
of the corresponding simplicial complement is $\o\Lambda^{*,*}(E_{\{i_1,i_2,i_3\}}\PP)$. The graded vector space
\[
\displaystyle{\bigoplus_{\sigma}H_{*,\sigma}(\o\Lambda^{*,\sigma}(E_{\{i_1, i_2\,i_3\}}\PP), d)}
\]
has Poincar\`{e} series $x^6(1+3x+3x^2+x^3)$.
\end{enumerate}

   Thus
\[
H^*(\ZZ_{K}({S^2}, {S^1}), \kk)=
 \bigoplus_{\omega\in K}\left(\bigoplus_{\tau}H_{q,\tau}(\o\Lambda^{*,\tau}(E_\omega\PP), d)\right)
\]
has Poincar\`{e} series
\begin{align*}
  & (1+3x^3+3x^6+x^9)\\
+ & 6x^2(1+x+2x^3+2x^4+x^6+x^7) \\
+ & 12x^4(1+2x+x^2+x^3+2x^4+x^5)\\
+ & 8x^6(1+3x+3x^2+x^3)\\
= & 1+6x^2+9x^3+12x^4+36x^5+35x^6+36x^7+54x^8+27x^9
\end{align*}
The total Betti number of $\ZZ_{K}({S^2}, {S^1})$ is $216$.


\begin{thebibliography}{00}

\bibitem{AP} C.~Allday and V.~Puppe, {\it Cohomological Methords in Transformation Groups,}
 Cambridge Studies in Advanced Mathematics, {\bf 32}, Cambridge University Press, 1993.
\bibitem{BB} A.~Bahri, M.Bendersky, F.~R.~Cohen and S.~Gitler, Decompositions of the polyhedral
 product functor with applications to moment-angle complexes and related spaces,
 {\it Proc. Nat. Acad. Sci. U. S. A.} {\bf 106} (2009), 12241-12244.
\bibitem{BBCG08} A.~Bahri, M.~Bendersky, F.~R.~Cohen and S.~Gitler,
{\it The polyhedral product functor: A methord of computation for moment-angle complexes,
arrangement and related spaces,}  preprint arXiv:0711.4689v2 [math.AT] 8 Dec 2008.
\bibitem{Baskakov} I.~Baskakov, {\it Cohomology of K-powers of spaces and the combinatorics
of simplicial divisions,} Russian Math. Surveys {\bf 57} (2002), no. 5, 989-990.
\bibitem{BP} V.~M.~Buchstaber and T.~E.~Panov, {\it Torus Actions and Their Applications in
 Topology and Combinatorics}, University Lecture Series, Vol. {\bf 24}, Amer. Math. Soc.
 Providence, RI, 2002.
\bibitem{BP04} V.~M.~Buchstaber and T.~E.~Panov, Combinatorics of simplicial cell complexes
 and tours action, {\it Proc. Steklov Inst. Math.} {\bf 247} (2004), 33-49.
\bibitem{CL} X.~Cao and Z.~L\"{u},  M\"{o}bius transform, moment-angle complexes and
 Halperin-Carlsson conjecture, {\it preprint.} arXiv:0908.3174v2 [math.CO] 12 Sep 2009.
\bibitem{DJ} M.~W.~Davis and T.~Januszkiewicz,  Convex polytopes, coxeter orbifolds and
 torus action, {\it Duck Math. J.} {\bf 62} (1991), 417-451.
\bibitem{DS} G.~Denham and A.~Suciu, Moment-angle complexes, monomial ideals and Massey products,
 {\it Pure Appl. Math. Q.} {\bf 3} (2007), 25-60.
\bibitem{Franz} M.~Franz, {The intergral cohomology of toric manifolds,} Proc. Steklov Inst. Math.
{\bf 252} (2006), 53-62. [Proceedings of the Keldysh Conference, Moscow 2004].
\bibitem{GM} R.~Goresky and R.~MacPherson, {\it Stratified Morse Theory,} Ergeb. Math. Grenzgeb.,
Vol. 14, Springer-Verlag, Berlin, 1988.
\bibitem{GTh} J.~Grbic and R.~Theriault, {\it Homotopy type of the complement of a coordinate
sunspace arrangement of codimension two,} Russian Math. Surveys {\bf 59} (2004), no. 3, 1207-1209.
\bibitem{Hilton} P.~J.~Hilton, U.~Stammbach, {\it A Course in Homological Algebra,}
 Berlin-Heideberg-New York: Springer 1971.
\bibitem{Hochster} M.~Hoster, {\it Cohen-Macaulay ring, combinatorics and simplecial complexes,}
in: Ring theory, II (Proc. Second Conf. Univ. Oklahoma, Norman, Okla., 1975), pp. 171-223.
\bibitem{Lopez} S.~Lopez de Medrano, {\it Topology of the intersection of quadrics in $\mathbb{R}^n$,}
in {\it Algebraic Topology} (Arcata Ca), Springer Verlag LNM {\bf 1370} (1989), Springer Verlag.
\bibitem{LP} Z.~Lu and T.~Panov, Moment-angle complexes from simplicial posets, {\it preprint}
arXiv:0912.2219v1 [math.AT] 11 Dec 2009.
\bibitem{Mil} E.~Miller and B.~Sturmfels, {\it Combinatorial Commutative Algebra},
 Graduate Texts in Math. {\bf 227}, Springer, 2005.
\bibitem{Panov} T.~E.~ Panov, Cohomology of face rings and tours actions, {\it London
 Math. Soc. Lect. Notes Ser.} {\bf 347} (2008), 165-201.
\bibitem{Porter} G.~Porter, {\it The homotopy groups of wedge of suspensions,} Amer. J. Math.,
{\bf 88} (1966), 655-663.
\bibitem{Stanley} R.~P.~Stanley {\it Combinatorics and Commutative Algebra} second edition, Progress in Math.
{\bf 41} Birkhauser, Boston, 1996.
\bibitem{U} Y.~Ustinovsky, Toral rank conjecture for moment-angle complexes, {\it preprint}
arXiv:0909.1053v2 [math.AT] 29 Sep 2009.
\bibitem{Vogt} R.~Vogt, {\it Homotopy limits and colimits,} Math. Z. {\bf 134} (1973), 11-52.
\end{thebibliography}
\end{document}